\documentclass[12pt]{article}

\usepackage{amsfonts,amssymb,amsbsy, amsmath, dsfont}
\usepackage[german,english]{babel}
\usepackage[latin1]{inputenc}
\usepackage{graphicx}
\usepackage{amsthm}
\usepackage{xcolor}
\usepackage{capt-of}

\usepackage[hyperindex,backref]{hyperref} 

\newcommand{\E}{{\mathbb E}}

\newcommand{\N}{{\mathbb N}}

\renewcommand{\P}{{\mathbb P}}
\newcommand{\Q}{{\mathbb Q}}
\newcommand{\R}{{\mathbb R}}

\newcommand{\Z}{{\mathbb Z}}

\newcommand{\cB}{{\cal B}}
\newcommand{\cC}{{\cal C}}

\newcommand{\ind}{\mathds{1}}

\newtheoremstyle{thm}
	{8pt}
	{8pt}
	{\itshape}
	{}
	{\bfseries }
	{}
	{\newline}
	{}%

\newtheoremstyle{namedthm}
	{8pt}
	{8pt}
	{\itshape}
	{}
	{\bfseries }
	{}
	{5pt}
	{}%
	
\newtheoremstyle{def}
	{8pt}
	{8pt}
	{}
	{}
	{\bfseries }
	{}
	{\newline}
	{}%

\newtheoremstyle{nameddef}
	{8pt}
	{8pt}
	{}
	{}
	{\bfseries }
	{}
	{5pt}
	{}%
	
\swapnumbers
\theoremstyle{thm}
\newtheorem{Thm}{Theorem}[section]

\theoremstyle{namedthm}

\theoremstyle{def}

\theoremstyle{nameddef}

\newlength{\XWidth}

\makeatletter
\def\moverlay{\mathpalette\mov@rlay}
\def\mov@rlay#1#2{\leavevmode\vtop{%
   \baselineskip\z@skip \lineskiplimit-\maxdimen
   \ialign{\hfil$\m@th#1##$\hfil\cr#2\crcr}}}
\newcommand{\charfusion}[3][\mathord]{
    #1{\ifx#1\mathop\vphantom{#2}\fi
        \mathpalette\mov@rlay{#2\cr#3}
      }
    \ifx#1\mathop\expandafter\displaylimits\fi}
\makeatother

\begin{document}

\author{Sebastian Ziesche\thanks{Karlsruhe Institute of Technology, sebastian.ziesche@kit.edu}
}

\title{Sharpness of the phase transition and lower bounds for the critical intensity in continuum percolation on $\R^d$}
\date{\today}
\maketitle

\begin{abstract}
\noindent We consider the Boolean model $Z$ on $\R^d$ with random compact grains, i.e. $Z := \bigcup_{i \in \N} (X_i + Z_i)$ where $\eta_t := \{X_1, X_2, \dots\}$ is a Poisson point process of intensity $t$ and $(Z_1, Z_2, \dots)$ is an i.i.d. sequence of compact grains (not necessarily balls). We will show, that the volume and diameter of the cluster of a typical grain in $Z$ have an exponential tail if the diameter of the typical grain is a.s. bounded by some constant. To achieve this we adapt the arguments of \cite{duminil2015newproof} and apply a new construction of the cluster of the typical grain together with arguments related to branching processes. \\
In the second part of the paper, we obtain new lower bounds for the Boolean model with deterministic grains. Some of these bounds are rigorous, while others are obtained via simulation. The simulated bounds are very close to the "true" values and come with confidence intervals.
\end{abstract}

\maketitle

\begin{flushleft}
\end{flushleft}

\section{Introduction}\label{sec:introduction}
Let $\cC^d$ be the set of nonempty compact subsets of $\R^d$ equipped with the usual Fell-topology and $\Q$ a measure on $\cC^d$. Let $\xi_t := \{(X_1, Z_1), (X_2, Z_2), \dots\}$ be a Poisson point process on $\R^d \times \cC^d$ with intensity measure $t \lambda^d \otimes \Q$, where $\lambda^d$ is the $d$-dimensional Lebesgue measure and $t \geq 0$. This corresponds to an independently marked Poisson process on $\R^d$ with intensity $t$ and mark distribution $\Q$. Without loss of generality we assume, that $\Q$ is concentrated on the grains that contain the origin. As a serious restriction of the model, we also have to assume, that there is a radius $R \in \R$ such that $\Q$ is concentrated on grains contained in the ball $B_R$ of radius $R$ centered at the origin (in the $2$-dimensional case, this assumption might be relaxed; see \cite{ahlberg2016sharpness}).

We write $\mathbf{N}(\R^d \times \cC^d)$ for the set of simple locally finite counting measures on $\R^d \times \cC^d$. We identify each element $\eta = \sum_{i = 1}^\tau \delta_{(x_i, K_i)} \in \mathbf{N}(\R^d \times \cC^d)$, $\tau \in \N \cup \{\infty \}$ with its support $\{(x_i, K_i) \mid i \in \N,\ i \leq \tau\}$. The Boolean model
\begin{align}\label{eq:def-BM}
  Z(\eta) := \bigcup_{i \in \N} (x_i + K_i)
\end{align}
is closely related to the graph $G(\eta) := (V, E)$ where $V := \{x_i \mid i \in \N\}$ and two distinct points $x_i, x_j$ are adjacent iff $(x_i + K_i) \cap (x_j + K_j) \neq \emptyset$, i.e. if the corresponding grains overlap. In this way the connected components of the Boolean model $Z(\xi_t)$ correspond to the connected components of $G(\xi_t)$ if $\Q$ is concentrated on the connected sets (this is true in most cases, but we will not need this assumption).


We want to study the connected component of a typical point of $\xi_t$ in $G(\xi_t)$ which is (due to the well known properties of the Palm distribution of Poisson processes, which can be found in \cite{schneider2008stochastic}) the same as studying the cluster of the origin in $G(\xi_t + \delta_{(0, Z_0)})$ where $Z_0 \sim \Q$ independently of $\xi_t$. Hence we define the \emph{zero cluster} $C_0 := C_0(t)$ as the set of all grains $(x,K) \in \xi_t$ where $x$ is connected to $0$ in $G(\xi_t + \delta_{(0, Z_0)})$. For $K \in \cC^d$ we write
\begin{align}\label{eq:def-radius}
  \rho(K) := \max_{x \in K} \|x\|_2
\end{align}
for the "radius" of $K$. We denote by $B_r(x)$ the euclidian ball centered in $x \in \R^d$ of radius $r$ and write $B_r := B_r(0)$.

The structure of the paper is as follows. In the second section we recall the tools needed to work with Poisson processes and discuss briefly the existence of a non-trivial phase transition in the model. In the third section we show the exponential decay of the tail of $\rho(C_0) := \max_{(x,K) \in C_0} \rho(x+K)$ and $|C_0|$ in the subcritical regime by constructing $C_0$ in a new way and using a comparison with a branching process. We will also show the mean-field lower bound $\P[|C_0(t)| = \infty ] \geq \frac{t - \tilde t_c}{t^2 \lambda^d(B_{2R})} $ for the percolation function in the supercritical regime. This is done by an adaptation of the arguments in \cite{duminil2015newproof}. In the fourth and fifth section we use these results to obtain new lower bounds for the critical intensity $t_c$ in the model where the grains are a.s. balls of radius one.

\section{Preliminaries}\label{sec:preliminaries}
To work with the Poisson process $\xi_t$ we need two well known tools. The first one is the Mecke-equation for Poisson processes which can be found in \cite[Thm 4.1]{last2016poissontobepublished}. Let $f:\mathbf{N}(\R^d \times \cC^d) \times (\R^d \times \cC^d) \to [0, \infty )$ be measurable, then
\begin{align}\label{eq:mecke}
  \E \sum_{(x, K) \in \xi_t} f(\xi_t, x,K) = t \int_{\cC^d} \int_{\R^d} \E[f(\xi_t + \delta_{(x,K)}, x, K)] \ dx \ \Q(dK).
\end{align}
Let $A \subset \mathbf{N}(\R^d \times \cC^d)$ be measurable with respect to the usual $\sigma$-algebra (see \cite{schneider2008stochastic}). The event $A$ is \emph{determined by a set $D \subset \R^d \times \cC^d$} if for all $\eta_1, \eta_2 \in \mathbf{N}(\R^d \times \cC^d)$ with $\eta_1 \cap D = \eta_2 \cap D$ we have
\begin{align}
  \eta_1 \in A \quad \Leftrightarrow \quad \eta_2 \in A.
\end{align}
The second important relation is a Margulis-Russo type formula for Poisson processes which can be found in \cite{last2014perturbation}. Let $D \subset \R^d \times \cC^d$ be such that $(\lambda^d \otimes \Q)(D) < \infty $ and let $A$ be an event that is determined by $D$, then
\begin{align}\label{eq:russo-formula}
  \frac{\partial \P[\xi_t \in A]}{\partial t} = \int_{\cC^d} \int_{\R^d} \E[\ind_A(\xi_t + \delta_{(x,K)}) - \ind_A(\xi_t)] \ dx \ \Q(dK).
\end{align}
A thorough treatment of the Poisson process can also be found in \cite{last2016poissontobepublished}.

We recall the definition of the critical intensity
\begin{align}
  t_c := \sup\{t \geq 0 \mid \P[|C_0(t)| = \infty ] = 0 \}
\end{align}
and want to point out, that under the assumption that $Z_0 \subset B_R$ a.s.\ it is easy to show, that $t_c > 0$. This is due to a simple coupling of the model with the model where $\Q(\{B_R\}) = 1$. However it is possible, that $t_c = \infty $ if the grains don't contain a small ball with positive probability (see \cite{hall1985continuum} for a more elaborate version of these statements). This would be the case if we had line segments as grains for example. But as there are also interesting grain distributions with lower dimensional grains, we don't want to exclude this case. The results are unaffected by that.

\section{A sharp phase transition}\label{sec:sharp_phase_transition}
For $D \subset \R^d$ we write $[D] \subset \cC^d$ for the measurable subset of grains, that intersect $D$, we write $\partial D$ for the set of grains that intersect $D$ as well as $D^c$ and we define $D^\circ := [D^c]^c$ the set of grains that are contained in $D$. For $D_1, D_2 \subset \R^d$ and $A \subset \cC^d$ measurable, the event "$D_1 \leftrightarrow D_2 \text{ in } \xi_t \cap A$" holds, iff there is a path $(X_{a_1}, \dots, X_{a_n})$ in $G(\xi_t)$ with $X_{a_1} + Z_{a_1} \cap D_1 \neq \emptyset $, $X_{a_n} + Z_{a_n} \cap D_2 \neq \emptyset $ and $X_{a_i} + Z_{a_i} \in A$ for all $i \in [n]$. With $\{D_1 \leftrightarrow \infty \text{ in } \xi_t \}$ we mean the event, that $D_1$ is intersected by an infinite cluster of $G(\xi_t)$ and remark, that this is equal to $\bigcap_{n \in \N} \{D_1 \leftrightarrow B_n^c \text{ in } \xi_t\}$.

The heart of the proof in \cite{duminil2015newproof} is the study of a functional $\varphi$ on subsets $S$ of $\Z^d$ containing the origin. The functional is equal to the expected number of open edges in the edge-boundary of $S$ that are connected to the origin in $S$. The proper counterpart of $\varphi$ in our model is defined for each $S \in \cB(\R^d)$ with $B_R \subset S$ by
\begin{align}
\begin{aligned}\label{eq:def_phi}
  \varphi_t(S) &:= \E \sum_{(x,K) \in \xi_t} \ind\{x+K \in \partial S,\ B_R \leftrightarrow x+K \text{ in } \xi_t \cap S^\circ\} \\
  &\ =t \int_{\cC^d} \int_{\R^d} \ind\{x+K \in \partial S\} \P[B_R \leftrightarrow x+K \text{ in }\xi_t \cap S^\circ]\ dx \ \Q(dK).
\end{aligned}
\end{align}
This is the expected number of grains that "cross the boundary of $S$" and are connected to $B_R$ by grains contained in $S$. We proceed as in \cite{duminil2015newproof} by defining a "new" critical intensity
\begin{align}
  \tilde t_c := \sup\{t \geq 0 \mid \exists S \in \cB(\R^d):\ B_R \subset S,\ \varphi_t(S) < 1\}.
\end{align}

\begin{Thm}\label{thm:sharp_phase_transition}
  We have $t_c = \tilde t_c$. Moreover
  \begin{enumerate}
    \item for $t < t_c$ there are constants $c_1, c_2 > 0$ depending on $t$ such that
        \begin{align}
          \P[\rho(C_0(t)) \geq r] & \leq e^{-c_1 r},\quad r \geq 0\\
          \P[|C_0(t)| \geq r] & \leq e^{-c_2 r},\quad r \geq 0
        \end{align}
    \item for $t = t_c$
        \begin{align}
          \E |C_0(t)| = \infty ,
        \end{align}
    \item for $t > t_c$
        \begin{align}
          \P[B_R \leftrightarrow \infty \text{ in } \xi_t ] & \geq \frac{t - \tilde t_c}{t},\\
          \P[|C_0(t)| = \infty ] & \geq \frac{t - \tilde t_c}{t^2 \lambda^d(B_{2R})}.
        \end{align}
  \end{enumerate}
\end{Thm}
For the Boolean model where $\Q$ is concentrated on $B_1$, Penrose gave a simple lower bound for the critical intensity in \cite{penrose1996continuum}. This was done by modifying the construction of the cluster $C_0$ in a monotone way.

To be more precise, Penrose showed, that the cluster of the Ball at the origin can be constructed by the following algorithm. Let $C = \emptyset $ and $D = \{(0,B_1)\}$. Do the following until $D$ is empty: take the first (via first in, first out procedure) element $(x,K)$ of $D$, create a Poisson process on the set of balls that intersect $(x,K)$ but none of the grains in $C$. Add these newly created grains to $D$ and move $(x,K)$ from $D$ to $C$. When the algorithm stops, $C$ has the same distribution as $C_0$.

To obtain the lower bound, Penrose removed the condition, that the newly created grains should not intersect $C$. This way more grains are added to $D$ and hence $C$ gets larger. However the algorithm also becomes simpler, as it don't cares about the past anymore and the number of grains in $C$ behaves like the total progeny of a Galton-Watson process where the number of offsprings is Poisson-distributed with parameter $t \lambda^d(B_2)$. If the intensity of the Poisson process is such that the expected number of offsprings is equal to one ($t = \lambda^d(B_2)^{-1}$), the corresponding Galton-Watson process almost surely dies out. Furthermore it is easy to show, that for any smaller intensity, the radius and volume of $C_0$ have an exponential tail (see the proof of \ref{thm:sharp_phase_transition}).

We will now refine this approach such that it is applicable to the whole subcritical regime and random grains.

\emph{Proof of Theorem \ref{thm:sharp_phase_transition}:} If we prove items 1.-3. for $\tilde t_c$ instead of $t_c$ the theorem follows immediately. For the first part, let $t < \tilde t_c$ and let $S \in \cB(\R^d)$ be such that $\varphi_t(S) < 1$ and $B_R \subset S$.

Before going into technical details, we want to give an informal description how we construct the cluster $C_0$ in our modification of the algorithm mentioned above. After constructing $Z_0$ we explore the cluster of $Z_0$ in $S^\circ$ (in the same way a in the original algorithm). Then we construct all grains contained in $\partial S$, that are connected to $Z_0$ in $S^\circ$. For each of these grains $(x,K)$ we repeat this procedure with $S$ being replaced by $S + x$ except, that we realize our Poisson process only once everywhere. We repeat this until there are no more grains that intersect a translated $S$ and its complement.

\begin{center}
  \includegraphics{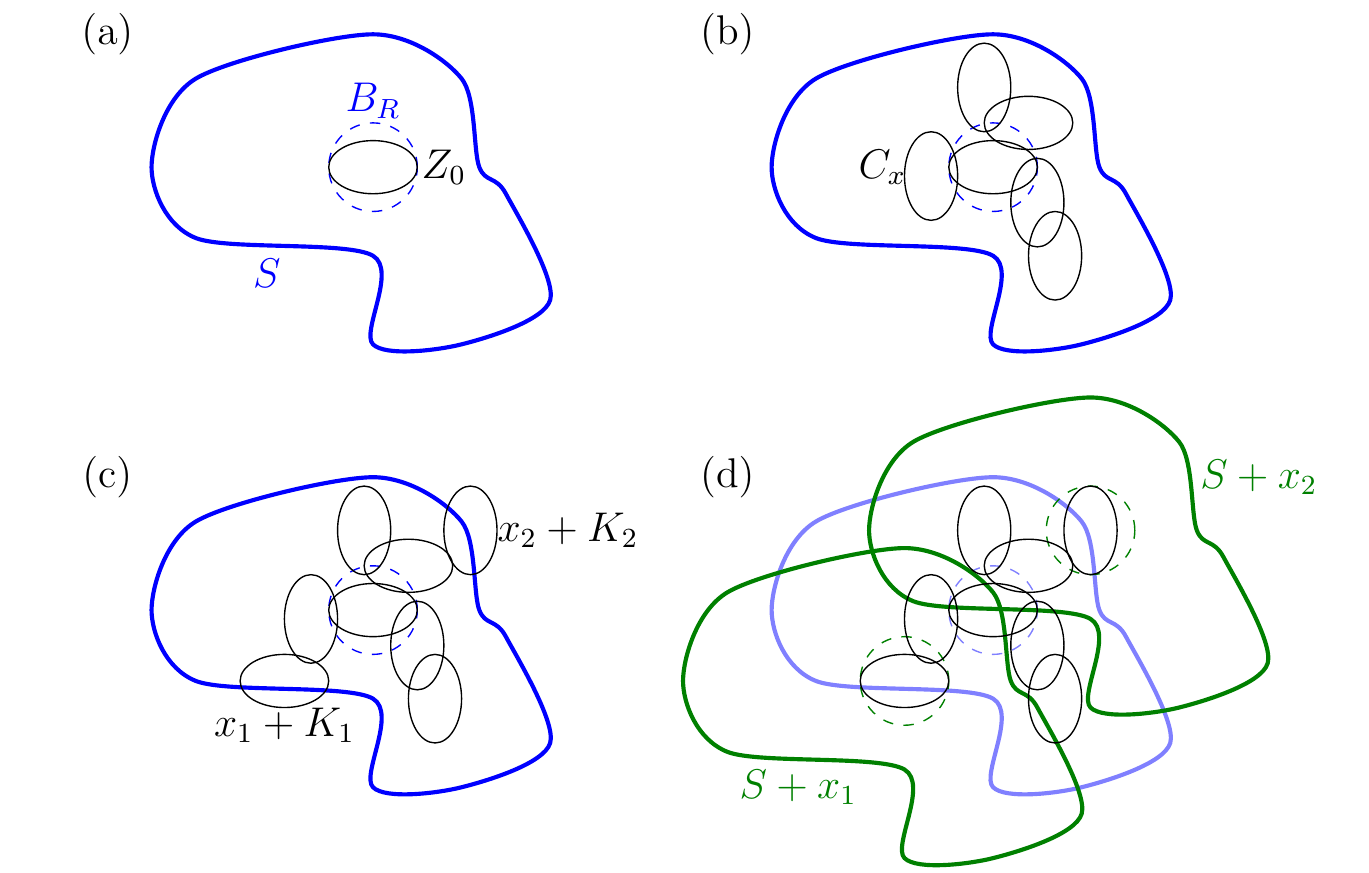}
  \captionof{figure}{The first steps of the new algorithm to construct $C_0$ where $\Q$ has probability mass $.5$ on each of two ellipses. In (a) the grain $Z_0$ at the origin is constructed. In (b) the cluster of $(0,Z_0)$ in $S^\circ$ is constructed. In (c) the grains in $\partial S$ that intersect $C_x$ are generated. (d) shows in green the translated versions of $S$ where the steps (b) and (c) are repeated with the only restriction, that newly generated grains must not intersect any previously generated grains other than the ones with center at $x_1$ or $x_2$.}
  \label{fig:1}
\end{center}

To formalize this, we use the following algorithm. Let $C = \emptyset $ and $D = \{(0,Z_0)\}$. Do the following until $D$ is empty: Take the first (via first in, first out procedure) element $(x,K)$ of $D$. Construct the cluster $C_{x}$ of $K+x$ with the grains that are contained in $(S+x) \setminus \bigcup_{(y,L) \in C} (y+L)$, in the same step by step manner as in Penrose' algorithm. Construct a Poisson process $\eta_x$ with intensity measure $t \lambda^d \otimes \Q$ restricted to the set $\partial (S + x) \cap [C_x] \cap [\bigcup_{(y,L) \in C} (y+L)]^c$. Add $C_x$ to $C$, remove $(x,K)$ from $D$ and add all grains in $\eta_x$ to $D$. This way $C$ has again the same distribution as $C_0$ when the algorithm stops.

Now we modify this algorithm in such a way that the resulting $C$ may not become smaller. First, we replace every grain $(x,K)$ that is put into $D$ by the grain $(x,B_R)$ (which contains $(x,K)$). Second, we construct the cluster $C_x$ not only with the grains contained in $(S+x) \setminus \bigcup_{(y,L) \in C} (y+L)$, but with the grains contained in $S + x$. This modification is in the same spirit as in the Penrose algorithm where the algorithm doesn't care anymore, what happened in the past.

As $\varphi_t(S)$, the expected number of grains in $\xi_t \cap \partial S$, that are connected to $B_R$ is less than one, it follows, that the modified algorithm will almost surely terminate. The cluster will then be related to the following Galton-Watson process $W$. The total progeny of $W$ consists of all the grains $(x,K)$ that were in $D$ at some point in time. A grain $(x,K)$ is a child of a grain $(y,L)$ if $(x,K)$ was created as one of the grains in $\partial(S+y)$ that was connected to $y+L$ in $S+y$ and hence added to $D$. It follows, that the offspring distribution of $W$ process has an expected value of $\varphi_t(S) < 1$ and that it is stochastically dominated by the total number of grains in a Poisson process of intensity $t \lambda^d \otimes \Q$ that lie in $\partial S$, which is a Poisson random variable.

Furthermore, the algorithm implies, that the radius of $C$ is less than or equal to the number of generations $W$ survived, times $\rho(S) + R$, as with each generation, the grains are only generated at a distance of at most $\rho(S) + R$ of the ancestor grain in the previous generation. By a similar argument, we have that the total number of grains in $C$ is stochastically dominated by the sum of $X$ independent Poi$((t \lambda^d \otimes \Q)([S]))$-distributed random variables, where $X$ is the total progeny of $W$.

Hence it remains to show, that the extinction time and the total progeny of $W$ have an exponential tail. It is however well known (see the end of the first chapter in \cite{athreya1972branching}) that the time to extinction has an exponential tail if the expected number of offsprings is less than one. To obtain the exponential tail of the total progeny, it is useful to observe (see \cite{dwass1969total}), that the total progeny of a Galton-Watson process with offspring distribution $\mu$ is equal to the first time the random walk, where the increments $I$ are distributed such that $I+1 \sim \mu$, hits $-1$ if it started at $0$. The increments of this random walk have expected value $\varphi_t(S) - 1 < 0$ and are stochastically dominated by a Poisson random variable. Hence the basic theorems of large deviation (e.g. \cite[Theorem 27.3.]{kallenberg2002foundations}), give the exponential decay of the total progeny.

To prove 2. we observe, that due to \cite[Theorem 3.1]{last2014perturbation} $\varphi_t(S)$ is an analytic and hence continuous function in $t$ for a fixed $S$. It follows, that the set of parameters $t$ where there is a set $S \in \cB(\R^d)$ such that $B_R \subset S$ and $\varphi_t(S) < 1$ is open in the interval $[0,1]$. We deduce, that for any $S \in \cB(\R^d)$ such that $B_R \subset S$ we have $\varphi_{\tilde t_c}(S) \geq 1$ and hence
\begin{align*}
  \E |C_0(\tilde t_c)| & = \E \sum_{(x,K) \in \xi_{\tilde t_c}} \ind\{Z_0 \leftrightarrow x + K \text{ in } \xi_{\tilde t_c}\} \\
  & \geq \sum_{n = 1}^\infty  \E \sum_{(x,K) \in \xi_{\tilde t_c}} \ind\{(x,K) \in \partial B_{3nR}, Z_0 \leftrightarrow x + K \text{ in } \xi_{\tilde t_c}\} \\
  & \geq \sum_{n = 1}^\infty  \varphi_{\tilde t_c}(B_{3nR}) \\
  & = \infty .
\end{align*}

We prove assertion 3. with the help of equation \eqref{eq:russo-formula} applied to the event $\{B_R \leftrightarrow B_r^c \text{ in } \xi_t \cap [B_r]\}$ for $r > R$. We have that
\begin{align*}
  \frac{\partial \P[B_R \leftrightarrow B_r^c \text{ in } \xi_t \cap [B_r]]}{\partial t}
  & = \int_{\cC^d} \int_{\R^d} \E[\ind\{B_R \leftrightarrow B_r^c \text{ in } \xi_t + \delta_{(x,K)}\cap [B_r] \\
  & \hspace{2.5cm}\text{ but not in } \xi_t \cap [B_r]\}] \ dx\ \Q(dK).
\end{align*}
Let $C(B_r^c)$ be the union of all grains, that are connected to $B_r^c$ in $\xi_t \cap [B_r]$, then
it is easy to verify, that $\{B_R \leftrightarrow B_r^c \text{ in } \xi_t + \delta_{(x,K)}\cap [B_r] \text{ but not in } \xi_t \cap [B_r]\}$ holds if and only if $\{B_R \leftrightarrow x+K \text{ in } \xi_t \cap [B_r^c \cup C(B_r^c)]^c,\ B_R \notin [C(B_r^c)],\ x+K \in \partial (B_r^c \cup C(B_r^c)) \}$ holds. Hence
\begin{align*}
  & \frac{\partial \P[B_R \leftrightarrow B_r^c \text{ in } \xi_t \cap [B_r]]}{\partial t} \\
  & = \int_{\cC^d} \int_{\R^d} \E[\ind\{B_R \leftrightarrow x+K \text{ in } \xi_t \cap [B_r^c \cup C(B_r^c)]^c,\\
  &\hspace{3cm} B_R \notin [C(B_r^c)],\ x+K \in \partial (B_r^c \cup C(B_r^c)) \}] \ dx\ \Q(dK).
\end{align*}
We condition on the shape of $C(B_r^c)$ (denoting its distribution by $\P_{C(B_r^c)}$) and remark, that conditioned on the event $C(B_r^c) = A$ the grains in $\xi_t \cap [B_r^c \cup C(B_r^c)]^c$ are distributed like $\tilde \xi_t \cap [B_r^c \cup A]^c$, where $\tilde \xi_t$ is a Poisson process independent of $\xi_t$ but with the same intensity measure. Hence
\begin{align*}
  & \frac{\partial \P[B_R \leftrightarrow B_r^c \text{ in } \xi_t \cap [B_r]]}{\partial t} \\
  & = \int_{\cC^d} \int_{\R^d} \int \E[\ind\{B_R \leftrightarrow x+K \text{ in } \tilde \xi_t \cap [B_r^c \cup A]^c,\ B_R \notin [A],\\
  &\hspace{3cm} x+K \in \partial (B_r^c \cup A) \} \mid C(B_r^c) = A] \ \P_{C(B_r^c)}(dA)\ dx\ \Q(dK)\\
  & = \int \int_{\cC^d} \int_{\R^d} \P[B_R \leftrightarrow x+K \text{ in } \xi_t \cap [B_r^c \cup A]^c]\ \ind\{x+K \in \partial (B_r^c \cup A)\}\\
  &\hspace{3cm} \ dx\ \Q(dK)\ \ind\{B_R \notin [A]\}\ \P_{C(B_r^c)}(dA)\\
  & = \frac{1}{t} \int \varphi_t((B_r^c \cup A)^c)\ \ind\{B_R \notin [A]\}\ \P_{C(B_r^c)}(dA).
\end{align*}
For any $t \geq \tilde t_c$ this yields, that
\begin{align}
  \frac{\partial \P[B_R \leftrightarrow B_r^c \text{ in } \xi_t \cap [B_r]]}{\partial t} & \geq \frac{1 - \P[B_R \leftrightarrow B_r^c \text{ in } \xi_t \cap [B_r]]}{t}.
\end{align}
Dividing the inequality by $1 - \P[B_R \leftrightarrow B_r^c \text{ in } \xi_t \cap [B_r]]$, integrating it from $\tilde t_c$ to some $t > \tilde t_c$ and using the trivial inequality $\P[B_R \leftrightarrow B_r^c \text{ in } \xi_{\tilde t_c}\cap [B_r]] \geq 0$ we obtain for $t \geq \tilde t_c$ that
\begin{align}\label{eq:bound_for_theta_t(r)}
  \P[B_R \leftrightarrow B_r^c \text{ in } \xi_t \cap [B_r]] \geq \frac{t - \tilde t_c}{t}.
\end{align}
Taking the limit $r \to \infty $ yields, that
\begin{align*}
  \P[B_R \leftrightarrow \infty \text{ in } \xi_t ] \geq \frac{t - \tilde t_c}{t}.
\end{align*}
The last step is to relate the Palm probability $\P[|C_0(t)| = \infty ]$ to the probability, that $B_R$ intersects the infinite cluster. We use the well known formula for Palm probabilities (see \cite[Theorem 3.3.2 and Theorem 3.5.3]{schneider2008stochastic}), that implies
\begin{align*}
  \P[|C_0(t)| = \infty ] & = \frac{1}{t \lambda^d(B_{2R})} \E \sum_{(x,K) \in \xi_t} \ind\{x \in B_{2R},\ x+K \leftrightarrow \infty \} \\
  & \geq \frac{1}{t \lambda^d(B_{2R})} \P[\exists (x,K) \in \xi_t:\ x \in B_{2R},\ x+K \leftrightarrow \infty] \\
  & \geq \frac{1}{t \lambda^d(B_{2R})} \P[B_R \leftrightarrow \infty \text{ in } \xi_t ] \\
  & \geq \frac{t - \tilde t_c}{t^2 \lambda^d(B_{2R})}.
\end{align*}
\qed

\section{Rigorous lower bounds for the critical intensity}\label{sec:rigorous_lower_bounds}
Theorem \ref{thm:sharp_phase_transition} opens up several ways to obtain lower bounds for the critical intensity. We get the first one, by choosing a certain $S$ and calculating or estimating $\varphi_t(S)$ depending on $t$. If we are able to choose $t$ such that $\varphi_t(S) < 1$, the chosen $t$ is a lower bound for $t_c$. In the special case, where $\Q$ is concentrated on $B_1$ we retrieve the Penrose bound, by setting $S := B_1$. Then $\varphi_t(S) = \E[|\xi_t \cap [B_1]|] = t \lambda^d(B_2)$ and hence $t_c \geq \lambda^d(B_2)^{-1}$.

But this exact lower bound can be sharpened by choosing $S := B_3$. In this case we have $\partial S = \{B_1(x) \mid 2 < \|x\|_2 \leq 4\}$, $S^\circ = \{B_1(x) \mid \|x\|_2 \leq 2\}$ and $R = 1$. Moreover any ball $B_1(y) \in S^\circ$ intersects $B_1$ and hence for any $x$ such that $B_1(x) \in \partial S$,
\begin{align*}
  \P[B_1 \leftrightarrow B_1(x) \text{ in } \xi_t \cap S^\circ ] & = \P[\xi_t([B_1] \cap [B_1(x)]) \geq 1]
  \\
  & = 1 - e^{-t \lambda^d(B_2 \cap B_2(x))}.
\end{align*}
It follows from definition \eqref{eq:def_phi}, that
\begin{align*}
  \varphi_t(S) & = t \int_{\R^d} \ind\{B_1(x) \in \partial S\} \P[B_1 \leftrightarrow B_1(x) \text{ in } \xi_t \cap S^\circ ]\ dx \\
  & = t \lambda^d(B_1) \int_{2}^4 r^{d-1} (1 - e^{-t \lambda^d(B_2 \cap B_2(r \mathbf{e}_1))})\ dr
\end{align*}
where $\mathbf{e}_1$ is an arbitrary vector of unit length. The integral may be evaluated numerically to almost arbitrary precision and hence it is easy to find a good approximation of the value $t$ where $\varphi_t(B_3) = 1$. The corresponding lower bounds are listed in the following table.
\begin{table}[t]
  \centering
  \begin{tabular}{lll}
    \hline\hline \noalign{\vskip 2mm}
    d & bound via $\varphi_t(B_3)$ & Penrose bound \\ \noalign{\vskip 2mm} \hline \noalign{\vskip 2mm}
    2 & 0.135802\dots & 0.0795774\dots\\
    3 & 0.0433691\dots & 0.0298415\dots\\
    4 & 0.0167131\dots & 0.0126651\dots\\
    5 & 0.00734445\dots & 0.00593678\dots\\
    6 & 0.00357261\dots & 0.00302358\dots\\
    7 & 0.00188850\dots & 0.00165352\dots\\
    8 & 0.00107117\dots & 0.000962435\dots\\
    9 & 0.000645942\dots & 0.000592123\dots\\
    10 & 0.000411202\dots & 0.000382941\dots\\
    11 & 0.000274803\dots & 0.000259158\dots \\ \noalign{\vskip 2mm} \hline \hline
  \end{tabular}
  \caption{Lower bounds for the critical intensity $t_c$ in dimension 2 to 11 when $\Q$ is concentrated on $B_1$.}
  \label{fig:exponential_decay_radius:phi_B_3_lower_bound}
\end{table}
These two approaches obviously work with other fixed grain shapes too, but the numerical calculations might become significantly more involved.

When working on lower bounds for the critical intensity, it has to be mentioned, that the first rigorous lower bound was given by Hall in \cite{hall1985continuum} for dimension $d = 2$. To our knowledge the approach has never been applied to higher dimensions and it was stated in the book by Meester and Roy \cite{meester1996continuum} that it was untractable there. We want to take a short moment to show that this is not the case and even gives bounds, that are better than the ones in Table \ref{fig:exponential_decay_radius:phi_B_3_lower_bound}.

The idea of Hall may be seen today as a refinement of the Penrose bound. Hall also modified the algorithm from Section \ref{sec:sharp_phase_transition} to construct $C_0$ by creating more Poisson points. In the algorithm the children of each grain $(x,K)$ are determined by the Poisson process of grains that intersect $x + K$ but not any previous grain in the cluster $C$. Halls modification was, to take all grains as children of $(x,K)$ that intersect $(x,K)$ but not the ancestor of $(x,K)$. In the case where all grains are a.s. equal to $B_1$ the number of children of the grain $B_1(x)$ only depends on the distance of $x$ to the center $y$ of its ancestor grain $B_1(y)$. In this way the cluster $C$ generated by the algorithm can be compared with a multitype branching process, where the type of each grain is the distance to the center of its ancestor grain. Hall showed, that the expected size of $C$ was given by
\begin{align}
  1 + \sum_{n = 1}^\infty t^n T^n( \ind_{(0,2)})(1)
\end{align}
where $T$ is an operator from the set of continuous functions on the interval $(0,2)$ onto itself, that is defined by
\begin{align}
  T(f)(x) := \int_{0}^2 f(y) g(y,x)\ dx.
\end{align}
The function $g(y,x)$ is given by the $d-1$ dimensional Hausdorff measure of the set $\{z \in \R^d \mid \|z\|_2 > 2,\ \|z - x \mathbf{e}_1\|_2 = y\}$. Hence to calculate $g$ in higher dimensions, we have to determine the surface area of a spherical cap in higher dimensions. The formulas for this can be found in the literature (see \cite{li2011concise}) and we obtain
\begin{align*}
  g(y, x) =
  \begin{cases}
    (d-1) \lambda^d(B_1) y^{d-1} \int_{0}^{\arccos\left( \frac{4- x^2 - y^2}{2 x y}\right)} \sin^{d-2}(\varphi) \ d\varphi, & y \in (2 - x, 2) \\
    0, & y \in (0, 2-x)
  \end{cases}.
\end{align*}
Hall concludes that the expected size of $C$ is finite, if the largest eigenvalue of $T$ is less than $1/t$. For fixed $d$ the integral defining $g$ can be solved analytically. Afterwards the largest eigenvalue of $T$ can be found numerically with very high precision. The results have been collected in Table \ref{tab:alle_schranken1}.

\section{Highly probable lower bounds for the critical intensity}
The second way to obtain lower bounds from Theorem \ref{thm:sharp_phase_transition} is, to use the mean-field lower bound. It follows from \eqref{eq:bound_for_theta_t(r)} by an easy calculation, that for any $t \in [0,\infty )$ and $r > 1$
\begin{align}
  t_c \geq t(1 - \theta_t(r)).
\end{align}
where
\begin{align}
  \theta_t(r) := \P[B_1 \leftrightarrow B_r^c \text{ in } \xi_t \cap [B_r]].
\end{align}
Hence to obtain a lower bound for $t_c$, it suffices to choose an arbitrary $t$ and $r > 1$ and estimate $\theta_t(r)$ by simulation. The event $\{B_1 \leftrightarrow B_r^c \text{ in } \xi_t \cap [B_r]\}$ can be simulated exactly and hence we may compute a rigorous one-sided confidence interval for the true value of $\theta_t(r)$ and hence a confidence interval for the lower bound of $t_c$.
\begin{table}[t]
  \centering
  \begin{tabular}{cllllll}
    \hline\hline \noalign{\vskip 2mm}
    $d$ & $r$ & $t$ & runs & success & 99\% CI for $\theta_t(r)$ & lower bound \\ \noalign{\vskip 2mm} \hline \noalign{\vskip 2mm}
    2 & 16000 & 0.357 & 10000 & 0 & 0.00063692 & 0.356772 \\
    3 & 2000 & 0.0814 & 10000 & 0 & 0.00063692 & 0.0813481\\
    4 & 500 & 0.0261 & 10000 & 10 & 0.002119993 & 0.0260445\\
    5 & 500 & 0.0101 & 10000 & 0 & 0.00063692 & 0.0100935\\
    6 & 200 & 0.00456 & 10000 & 1 & 0.000813077 & 0.00455628\\
    7 & 200 & 0.00228 & 10000 & 18 & 0.003154537 & 0.00227278\\
    8 & 150 & 0.00124 & 10000 & 21 & 0.003529665 & 0.00123560\\
    9 & 150 & 0.000725 & 10000 & 6 & 0.001571485 & 0.000723859\\
    10 & 120 & 0.000450 & 10000 & 4 & 0.001282615 & 0.000449422\\
    11 & 120 & 0.0002955 & 10000 & 8 & 0.001849554 & 0.000294952\\  \noalign{\vskip 2mm}\hline \hline
  \end{tabular}
  \caption{Simulation results for lower bounds of the critical intensity $t_c$ in dimension 2 to 11 in the Boolean model where $\Q$ is concentrated on $B_1$.}
  \label{tab:simulated_lower_bounds}
\end{table}

\begin{table}[t]
  \centering
  \begin{tabular}{llll}
    \hline\hline \noalign{\vskip 2mm}
    d & Penrose & via $\varphi_t(B_3)$  & Hall \\ \noalign{\vskip 2mm} \hline \noalign{\vskip 2mm}
    2 & 0.0795774 & 0.135802 & 0.174746 \\
    3 & 0.0298415 & 0.0433691 & 0.0534187 \\
    4 & 0.0126651 & 0.0167131 & 0.0198296 \\
    5 & 0.00593678 & 0.00734445 & 0.00845546 \\
    6 & 0.00302358 & 0.00357261 & 0.00401478 \\
    7 & 0.00165352 & 0.00188850 & 0.00208114 \\
    8 & 0.000962436 & 0.00107117 & 0.00116176 \\
    9 & 0.000592124 & 0.000645943 & 0.000691455 \\
    10 & 0.000382941 & 0.000411203 & 0.000435437 \\
    11 & 0.000259158 & 0.000274804 & 0.000288394 \\  \noalign{\vskip 2mm} \hline \hline
  \end{tabular}
  \caption{Rigorous lower bounds for $t_c$ from Section \ref{sec:rigorous_lower_bounds} for the Boolean model where $\Q$ is concentrated on $\{B_1\}$.}
  \label{tab:alle_schranken1}
\end{table}

\begin{table}[t]
  \centering
  \begin{tabular}{llllll}
    \hline\hline \noalign{\vskip 2mm}
    & Sim. $\theta_r(t)$, 99\% CI & Sim. \cite{torquato2012effect2} & Sim. \cite{torquato2012effect2} \\
    d & for lower bound & lower bound & upper bound \\ \noalign{\vskip 2mm} \hline \noalign{\vskip 2mm}
    2 & 0.356772 & 0.359076 & 0.359085\\
    3 & 0.0813481 & 0.081854 & 0.081858\\
    4 & 0.0260445 & 0.02632 & 0.02642\\
    5 & 0.0100935 & 0.01032 & 0.01034\\
    6 & 0.00455628 & 0.004516 & 0.004526\\
    7 & 0.00226708 & 0.002218 & 0.002272\\
    8 & 0.00123560 & 0.001206 & 0.001208\\
    9 & 0.000722539 & 0.0007121 & 0.0007133\\
    10 & 0.000449422 & 0.0004450 & 0.0004462\\
    11 & 0.000294952 & 0.0002933 & 0.0002935 \\  \noalign{\vskip 2mm} \hline \hline
  \end{tabular}
  \caption{Simulated lower bounds for $t_c$ for the Boolean model where $\Q$ is concentrated on $\{B_1\}$ compared to the best known values in the literature.}
  \label{tab:alle_schranken2}
\end{table}
We also want to point out, that for any $t < t_c$ the limit $\lim_{r\to \infty } \theta_t(r) = 0$. Hence by investing enough computing time it is in principle possible to approximate the value of $t_c$ arbitrarily well.

We did our simulations in the following way. We fixed a dimension, chose $t$ slightly below the best known value for $t_c$ from the literature and picked $r$ such that our simulations could finish in reasonable time. We simulated 10000 times the cluster $C_0$ with the algorithm stated before the proof of Theorem \ref{thm:sharp_phase_transition} and counted the number of times it intersected $B_r^c$. In higher dimensions we also terminated the algorithm when the size of $C$ exceeded some large threshold. In this case we counted this run as if $C$ had intersected $B_r^c$ and hence had a conservative estimate. After 10000 runs, we computed the corresponding confidence interval for $\theta_t(r)$ with the \texttt{prop.test} method of the statistical language R and calculated the lower bound for $t_c$. Depending on how fast this was done and on how many times the boundary was reached, we increased $t$ and started another 10000 runs. This lead to the results in Table \ref{tab:simulated_lower_bounds}.

We want to briefly discuss the chosen parameters and results. It can be observed, that the precision never exceeds three significant digits. This is due to the fact, that we only get a high precision, if the confidence interval is small. The size of the confidence interval however depends on the number of runs and successes. It turned out in practice, that choosing 10000 runs and a $t$ such that not more than about 20 runs succeed, gave the best tradeoff between time and precision. The few runs, where the cluster actually reaches $B_r^c$ are extremely time consuming, hence it is more efficient, to chose $t$ slightly below the expected "true" $t_c$. Nevertheless it can be seen in table \ref{tab:alle_schranken2}, that for high dimensions our lower bounds exceed the upper error bound (upper end of the 1-$\sigma$ band) of best simulation results in the literature (see \cite{torquato2012effect2}).

A few words concerning the implementation of the algorithm. It is very useful to save the approximate position of the grains in the cluster to have a faster access when comparing, if the current grain intersects $C$. Due to the fact that only a tiny fraction of the space is covered by the cluster, when $t$ is close to critical, we preferred a hashmap over an array for this task. Another important issue is the ball-picking method, i.e. the method to generate a random vector in $B_1$. This can either be done by generating points in $[-1, 1]^d$ and throwing away the points that don't lie in $B_1$ or it can be done by the formula proposed in \cite{barthe2005probabilistic}. The first approach is faster than the second one in low dimensions. We found that the second method was faster for $d \geq 7$.

\bibliography{C:/Arbeit/LaTeX/myBib}

\begin{thebibliography}{BGMN05}

\bibitem[AN72]{athreya1972branching}
K.~B. Athreya and P.~E. Ney.
\newblock {\em Branching Processes}, volume 196.
\newblock Springer Science \& Business Media, 1972.

\bibitem[ATT16]{ahlberg2016sharpness}
Daniel Ahlberg, Vincent Tassion, and Augusto Teixeira.
\newblock Sharpness of the phase transition for continuum percolation in
  $\mathbb{R}^{2}$.
\newblock {\em arXiv preprint arXiv:1605.05926}, 2016.

\bibitem[BGMN05]{barthe2005probabilistic}
F.~Barthe, O.~Guédon, S.~Mendelson, and A.~Naor.
\newblock A probabilistic approach to the geometry of the lp-ball.
\newblock {\em Ann. Probab.}, 33(2):480--513, 03 2005.

\bibitem[DCT15]{duminil2015newproof}
H.~Duminil-Copin and V.~Tassion.
\newblock A new proof of the sharpness of the phase transition for {B}ernoulli
  percolation on $\mathbb{Z}^{d}$.
\newblock {\em arXiv preprint arXiv:1502.03051}, 2015.

\bibitem[Dwa69]{dwass1969total}
M.~Dwass.
\newblock The total progeny in a branching process and a related random walk.
\newblock {\em Journal of Applied Probability}, 6(3):682--686, 1969.

\bibitem[Hal85]{hall1985continuum}
P.~Hall.
\newblock On continuum percolation.
\newblock {\em The Annals of Probability}, 13(4):1250--1266, 1985.

\bibitem[Kal02]{kallenberg2002foundations}
O.~Kallenberg.
\newblock {\em Foundations of Modern Probability}.
\newblock Springer Verlag, 2002.

\bibitem[Las14]{last2014perturbation}
G.~Last.
\newblock Perturbation analysis of {P}oisson processes.
\newblock {\em Bernoulli}, 20(2):486--513, 05 2014.

\bibitem[Li11]{li2011concise}
Shengqiao Li.
\newblock Concise formulas for the area and volume of a hyperspherical cap.
\newblock {\em Asian Journal of Mathematics and Statistics}, 4(1):66--70, 2011.

\bibitem[LP16]{last2016poissontobepublished}
G.~Last and M.D. Penrose.
\newblock {\em Lectures on the Poisson Process}.
\newblock 2016.
\newblock to be published by Cambridge University Press. Available at
  \url{http://www.math.kit.edu/stoch/~last/page/lehrbuch_poissonp/en}.

\bibitem[MR96]{meester1996continuum}
R.~Meester and R.~Roy.
\newblock {\em Continuum Percolation}.
\newblock Number 119 in Cambridge Tracts in Mathematics. Cambridge University
  Press, 1996.

\bibitem[Pen96]{penrose1996continuum}
M.~D. Penrose.
\newblock Continuum percolation and euclidean minimal spanning trees in high
  dimensions.
\newblock {\em The Annals of Applied Probability}, 6(2):528--544, 1996.

\bibitem[SW08]{schneider2008stochastic}
R.~Schneider and W.~Weil.
\newblock {\em Stochastic and Integral Geometry}.
\newblock Springer Verlag, 2008.

\bibitem[TJ12]{torquato2012effect2}
S.~Torquato and Y.~Jiao.
\newblock Effect of dimensionality on the continuum percolation of overlapping
  hyperspheres and hypercubes. ii. simulation results and analyses.
\newblock {\em The Journal of Chemical Physics}, 137(7), 2012.

\end{thebibliography}
\end{document}